\documentclass[aps,preprint,a4paper,floatfix]{revtex4}

\usepackage{graphicx,color}
\usepackage{amsmath,amssymb,amsfonts,amsthm,amscd,bm}
\usepackage{booktabs}
\usepackage{cancel}

\def\tilde{\widetilde}
\def\bar{\overline}
\def\hat{\widehat}
\def\*{\star}
\def\[{\left[}
\def\]{\right]}
\def\({\left(}      

\def\){\right)}

\def\frac#1#2{\dfrac{#1}{#2}}
\def\inv#1{\dfrac{1}{#1}}

\def\2pi{\hbox{$2\pi i$}}

\def\dsl{\raise.15ex\hbox{/}\kern-.57em\partial}
\def\Dsl{\,\raise.15ex\hbox{/}\mkern-.13.5mu D}
\def\sig{\sigma}

\def\2pi{\hbox{$2\pi i$}}

\def\dsl{\raise.15ex\hbox{/}\kern-.57em\partial}
\def\Dsl{\,\raise.15ex\hbox{/}\mkern-.13.5mu D}
\def\barray{\begin{eqnarray}}
\def\earray{\end{eqnarray}}
\def\beq{\begin{equation}}
\def\eeq{\end{equation}}

\def\AA{\leavevmode\setbox0=\hbox{h}
\dimen0=\ht0 \advance\dimen0 by-1ex\rlap{\raise.67\dimen0\hbox{\char'27}}A}

\makeatletter
\def\iddots{\mathinner{\mkern1mu\raise\p@
\vbox{\kern7\p@\hbox{.}}\mkern2mu
\raise4\p@\hbox{.}\mkern2mu\raise7\p@\hbox{.}\mkern1mu}}
\makeatother

\renewcommand*{\arraystretch}{1}
\renewcommand*{\arraycolsep}{5pt}
\makeatletter
\renewcommand*\env@matrix[1][\arraystretch]{%
  \edef\arraystretch{#1}%
  \hskip -\arraycolsep
  \let\@ifnextchar\new@ifnextchar
  \array{*\c@MaxMatrixCols c}}
\makeatother

\def\fref#1{FIG. \ref{#1}}
\def\tref#1{TABLE \ref{#1}}

\theoremstyle{plain}

\theoremstyle{remark}

\begin{document}

\title{
On the zeros of \boldmath{$L$}-functions
}
\author{
Guilherme Fran\c ca\footnote{guifranca@gmail.com}  
and  Andr\'e  LeClair\footnote{andre.leclair@gmail.com}
}
\affiliation{Cornell University, Physics Department, Ithaca, NY 14850} 

\begin{abstract}
We generalize our recent construction of the zeros of the  
Riemann $\zeta$-function to two infinite classes of $L$-functions,  
Dirichlet $L$-functions and those based on level one  modular forms. 
More specifically, we show that there are an infinite number of 
zeros on the critical line which are in 
one-to-one correspondence with the zeros of the cosine function,  
and thus enumerated by an integer $n$.
We obtain an exact equation
on the critical line that determines the $n$-th zero of these $L$-functions.  
We show that the counting formula on the critical line derived 
from such an equation agrees with the known counting formula on the 
entire critical strip.
We provide numerical evidence  supporting our statements, by computing
numerical solutions of this equation, yielding $L$-zeros to
high accuracy.    
We study in detail the $L$-function for the modular form based on  
the Ramanujan $\tau$-function, 
which is closely related to the bosonic string partition function.  
The same analysis for a more general class of $L$-functions 
is also considered.
\end{abstract}

\maketitle

\section{Introduction}

Dirichlet $L$-series are functions of a complex variable $z$ defined by
the series 
\beq
\label{Dseries}
L(z)  = \sum_{n=1}^\infty    \frac{a(n)}{n^z},  \qquad \Re(z) > 1
\eeq
where $a(n)$ is an arbitrary arithmetic function.     
In this paper we will consider two infinite classes of 
important $L$-functions,  
the Dirichlet $L$-functions where $a(n) = \chi (n)$ is a Dirichlet 
character,  and  $L$-functions associated with modular forms.   
The former have applications primarily in multiplicative number theory,
whereas the latter in additive number theory.  

The Dirichlet $L$-functions are generalizations of the 
Riemann $\zeta$-function,   
the latter being the simplest example \cite{Apostol}. 
They can be analytically continued to the entire 
complex plane.     The \emph{Generalized Riemann Hypothesis} (GRH)  is 
the conjecture that all non-trivial zeros 
of these functions are on the \emph{critical line},  i.e.  have real part 
equal to $\tfrac{1}{2}$.
Proving the validity of the GRH would have 
many  implications in number theory.     
Much less is known about the zeros of Dirichlet $L$-functions in 
comparison with the $\zeta$-function,  however let us mention a few works.   
Selberg \cite{Selberg1} obtained the analog of Riemann-von Mangoldt 
counting formula $N(T,\chi)$ for the number of zeros up to height $T$ within 
the entire critical strip $0 \le \Re(z) \le 1$. 
Based on this result, 
Fujii \cite{Fujii} gave an estimate for the number of zeros on the
critical strip with the ordinate between $[T+H, T]$. 
The distribution of
low lying zeros of $L$-functions near and at the critical line was examined
in \cite{Iwaniec}, assuming the GRH.
The statistics of the zeros,  i.e. the analog of the 
Montgomery-Odlyzko conjecture, were studied in 
\cite{Conrey,Hughes}.
It is also known that more than half of the non-trivial zeros are
on the critical line \cite{Conrey2}.
For a more detailed introduction to $L$-functions see \cite{Bombieri}.

Besides the Dirichlet $L$-functions, there are more general constructions
of $L$-functions based on arithmetic and geometric objects, 
like varieties over number fields and modular forms \cite{Iwaniec2, Sarnak}.
Some results for general $L$-functions are still conjectural. For instance,
it is not even clear if some $L$-functions can be analytically 
continued into a meromorphic function.   
We will only consider   the  additional $L$-functions
based on modular forms here.      
Thus the  $L$-functions considered in this paper have similar properties, 
namely, they posses an Euler product,  
can be analytic continued into the (upper half) complex plane, 
except for possible poles at $z=0$ and $z=1$, and satisfy a functional 
equation.

In our previous work \cite{GuiLeClair}, a new approach to the 
characterization of zeros of the $\zeta$-function
was developed,  building on the earlier work \cite{AL}.   
Enumerating the zeros on the critical line  as  
$\rho_n = \tfrac{1}{2}+ i y_n$,  with $n=1,2,\dotsc$, an exact 
transcendental  equation 
for the imaginary parts of 
the zeros $y_n$ was derived which depends only on $n$. An asymptotic 
version of these equations was first proposed in \cite{AL}.
From these equations for the zeros on the \emph{critical line}, one can derive 
the Riemann-van Mangoldt  and the exact 
Backlund counting formulae for the zeros on the \emph{entire strip}, 
thus this result strongly indicates that 
all non-trivial zeros of the $\zeta$-function are on the critical line. 
These transcendental equations can easily be 
solved numerically to very high precision, even for high zeros such as 
the billion'th,   with simple implementations such as are  available on
Mathematica.  Also, various approximate versions of the equation were 
presented.

In this work, we extend the results obtained in \cite{GuiLeClair}
to the two infinite  classes of $L$-functions mentioned above.  
In particular, we derive an exact equation for 
the $n$-th zero which only depends on $n$ and some 
very elementary 
properties of the mathematical object that the $L$-function is based on,
for instance, the Gauss sum of the Dirichlet character  
or the weight of its  modular form.
As for the $\zeta$-function,  the derivation presents an argument that  
all the zeros are on the critical line. 
We will not present our derivation in as much 
detail as we did for the  $\zeta$-function case \cite{GuiLeClair}, 
since our analysis follows precisely the same steps.

Our results are presented as follows. In section \ref{sec:character} 
we review some 
of the  main properties of Dirichlet characters and their  $L$-functions,  
all of which are well-known. In section \ref{sec:equation} we derive 
the exact equation for the $n$-th zero of Dirichlet $L$-functions
\eqref{exact}. One interesting new feature, in comparison with 
the $\zeta$-function, is that if the characters are complex numbers, 
then the zeros on the negative $y$-axis are not symmetrical to the ones 
lying on the positive $y$-axis. In section \ref{sec:approx},  we consider an 
approximation to the exact equation based on the leading order of the
generalized Riemann-Siegel $\vartheta_{k,a}$ function \eqref{RSgen} 
for large $y$. From this, one can derive counting formulas 
$N_0^{\pm}(T,\chi)$ for both  the positive and negative imaginary $y$-axes. 
An approximate formula for the zeros can be found explicitly  in terms of 
the Lambert $W$ function. In section \ref{sec:numerical}  we find 
numerical solutions to our equation \eqref{exact} to high accuracy,
considering two different examples of Dirichlet $L$-functions based on 
characters
with modulus $7$.  
In section \ref{sec:modular} we repeat the same steps but for 
$L$-functions based on modular forms, and we consider a specific example
based on the Ramanujan $\tau$-function, connected with
the Dedekind $\eta$-function and the bosonic string partition function.
 In appendix \ref{sec:general} we derive 
the main formulas for a more general, but not specific,  
class of $L$-functions.

\section{Dirichlet \boldmath{$L$}-functions}
\label{sec:dirichlet}

\subsection{Dirichlet characters and \boldmath{$L$}-functions}
\label{sec:character}

In this section we recall some of the main properties of Dirichlet 
characters and  $L$-functions based on it \cite{Apostol}.

Dirichlet $L$-series are defined as 
\beq
\label{Ldef}
L(z, \chi) = \sum_{n=1}^\infty \frac{\chi (n)}{n^z}, \qquad \Re (z) > 1
\eeq
where the arithmetic function $\chi(n)$ is a Dirichlet character.  
They  can all be analytically continued to the entire 
complex plane, and are then referred to as  Dirichlet $L$-functions. 

There are an infinite number of distinct Dirichlet 
characters which are primarily characterized by their modulus $k$,
which determines their periodicity. They can
be defined axiomatically,  which leads to specific properties,
some of which we now describe.
Consider a Dirichlet character $\chi$ mod $k$,  and let 
the symbol $(n,k)$ denote the greatest  common divisor of the
two integers $n$ and $k$. Then $\chi$  has the 
following properties:
\begin{enumerate}
\item $\chi(n+k) = \chi (n)$.
\item $\chi(1) = 1$ and  $\chi(0) = 0$.
\item $\chi(n  m) = \chi (n) \chi (m)$.
\item $\chi(n) = 0$ if  $(n, k) > 1$ 
and $\chi(n) \ne 0$ if $(n, k) = 1$.
\item \label{root1} If $(n, k) = 1$  then  $\chi(n)^{\varphi(k)} = 1$, where
$\varphi(k)$  is the Euler totient arithmetic function. This implies  
that $\chi(n)$ are roots of unity.
\item If $\chi$ is a Dirichlet character so is the complex conjugate $\chi^*$.
\end{enumerate}
For a given modulus $k$ there are $\varphi(k)$ distinct Dirichlet 
characters,  which essentially follows from property \ref{root1} above.
They can thus be labeled as $\chi_{k,j}$ where $j= 1, 2,\dotsc, \varphi(k)$ 
denotes an arbitrary ordering. 
If $k=1$ we have the \emph{trivial} character where
$\chi(n)=1$ for every $n$, and \eqref{Ldef} reduces to the Riemann 
$\zeta$-function.
The \emph{principal } character,  usually denoted $\chi_1$,   is defined as 
$\chi_1(n) = 1$ if $(n,k) = 1$ and zero otherwise. In the above notation 
the principal character is always $\chi_{k,1}$.

Characters can be classified as \emph{primitive} or \emph{non-primitive}.
Consider the Gauss sum
\beq\label{tau}
G(\chi) = \sum_{m=1}^{k}\chi(m)e^{2\pi i m / k}.
\eeq
If the character $\chi$ mod $k$ is primitive, 
then $|G(\chi)|^2 = k$. This is no longer valid for a non-primitive character.
Consider a non-primitive character $\bar{\chi}$ 
mod $\bar{k}$. Then it can be expressed in terms of
a primitive character of smaller modulus
as $\bar{\chi}(n) = \bar{\chi_1}(n) \chi(n)$, where $\bar{\chi_1}$ is
the principal character mod $\bar{k}$ and $\chi$ is a primitive
character mod $k < \bar{k}$, where $k$ is a divisor of $\bar{k}$. More
precisely, $k$ must be the \emph{conductor} of $\bar{\chi}$ 
(see \cite{Apostol} for further details).
In this case the two $L$-functions are related as
$L(z, \bar{\chi}) = L(z,\chi) \Pi_{p | \bar{k}}\(1 - \chi(p)/p^z\)$. 
Thus $L(z,\bar{\chi})$ has the same zeros as $L(z,\chi)$.
The principal character is only primitive when $k=1$, which yields the
$\zeta$-function. The simplest example of non-primitive characters
are all the principal ones for $k \ge 2$, whose zeros are the same as 
the $\zeta$-function. Let us consider another example with $\bar{k}=6$,
where $\varphi(6)=2$,  
namely $\bar{\chi}_{6,2}$, whose components are 
\footnote{Our enumeration convention for the $j$-index of $\chi_{k,j}$ 
is taken from Mathematica.} 
\beq
\begin{tabular}{@{}c|cccccc@{}}
$n$             & $1$ & $2$ & $3$ & $4$ & $5$ & $6$ \\
\midrule[0.3pt]
$\bar{\chi}_{6,2}(n)$ & $1$ & $0$ & $0$ & $0$ & $-1$ & $0$
\end{tabular}
\eeq
In this case, the only divisors are $2$ and $3$. Since $\chi_1$ mod $2$
is non-primitive, it is excluded. We are left with $k=3$ which is the
conductor of $\bar{\chi}_{6,2}$. Then we have  
two options; $\chi_{3,1}$ which is the non-primitive 
principal character mod $3$, thus excluded, 
and $\chi_{3,2}$ which is primitive. Its components are
\beq
\begin{tabular}{@{}c|ccc@{}}
$n$             & $1$ & $2$ & $3$ \\
\midrule[0.3pt]
$\chi_{3,2}(n)$ & $1$ & $-1$ & $0$ 
\end{tabular}
\eeq
Note that $|G(\chi_{6,2})|^2 = 3 \ne 6$ and $|G(\chi_{3,2})|^2 = 3$.
In fact one can check that 
$\bar{\chi}_{6,2}(n) = \bar{\chi}_{6,1}(n) \chi_{3,2}(n)$, 
where $\bar{\chi}_{6,1}$ is the principal character 
mod $\bar{k}=6$. Thus the zeros of $L(z,\bar{\chi}_{6,2})$ are the same as 
those of $L(z,\chi_{3,2})$. Therefore, it suffices to consider 
primitive characters, and we will henceforth do so.  

We will need the functional equation satisfied by $L(z,\chi)$.   
Let $\chi$ be a \emph{primitive} character. 
Define its \emph{order} $a$ such that 
\beq\label{order}
a \equiv \begin{cases}1 \qquad
\mbox{if $\chi(-1)= -1$ ~~(odd order)} \\
0 \qquad \mbox{if $\chi(-1) = 1$  ~~~~(even order)}
\end{cases}
\eeq
Let us define the meromorphic function 
\beq
\label{Lambda}
\Lambda(z,\chi) \equiv \( \dfrac{k}{\pi} \)^{\tfrac{z+a}{2}}   
\Gamma\(\dfrac{z+a}{2}\) L(z, \chi).   
\eeq
Then $\Lambda$  satisfies the well known functional equation \cite{Apostol}
\beq
\label{FELambda}
\Lambda (z, \chi) =  \dfrac{ i^{-a}  \, G(\chi)}{\sqrt{k}} 
\Lambda (1-z, \chi^*).
\eeq
The above equation is only valid for primitive characters.

\subsection{The exact  transcendental equation for the \boldmath{$n$}-th 
zero on the critical line}
\label{sec:equation}

In this section we derive an exact equation satisfied by zeros enumerated by 
an integer $n$. The analysis that leads to this
equation is the same as for the $\zeta$-function in \cite{AL,GuiLeClair},  
consequently we do not provide as detailed a derivation,  
since such details can be
surmised from \cite{GuiLeClair}.   

For a primitive character,  since $|G(\chi)| = \sqrt{k}$,  
the factor on the right hand side of \eqref{FELambda} is a phase. 
It is  thus possible to obtain a more symmetric form of 
\eqref{FELambda} through a new function defined as
\beq
\label{xi}
\xi(z, \chi) \equiv \dfrac{i^{a/2} \, k^{1/4} }{\sqrt{G\(\chi\)}} \, \, 
\Lambda (z, \chi). 
\eeq
It  then satisfies 
\beq
\label{FExi}
\xi(z, \chi)  = \xi^{*}(1-z , \chi ) \equiv   \(\xi(1-z^*, \chi)\)^*.
\eeq
Above, the function $\xi^*$  of $z$ is defined as the complex conjugation of
all coefficients that  define $\xi$,  
namely  $\chi$ and the $i^{a/2}$ factor, evaluated at a non-conjugated $z$.
   
%Note that unlike $\Lambda$ given by \eqref{Lambda}, which obeys 
%$(\Lambda(z, \chi))^* = \Lambda(z^*,\chi^*)$, for \eqref{xi} we now have
%\beq\label{xiconju}
%\(\xi(z,\chi)\)^* = (-1)^{a/2}\xi\(z^*,\chi^*\). 
%\eeq
Note that $(\Lambda(z, \chi))^* = \Lambda(z^*,\chi^*)$. Using the well
known result $G\(\chi^*\) = \chi(-1) \( G(\chi) \)^*$ we conclude that
\beq\label{xiconju}
\(\xi(z,\chi)\)^* = \xi\(z^*,\chi^*\).
\eeq
This implies that if the character is \emph{real}, then if $\rho$ is a 
zero of $\xi$  so is $\rho^*$, and one needs only 
consider $\rho$ with  
positive imaginary part. On the other hand if $\chi \neq \chi^*$, 
then the zeros with negative imaginary part are different than $\rho^*$.
For the trivial character where $k=1$ and $a=0$, implying 
$\chi(n) = 1$ for any $n$,  then $L(z,\chi)$  reduces to the
Riemann $\zeta$-function and  \eqref{FExi} yields the well
known functional equation $\xi(z)=\xi(1-z)$ with 
$\xi(z) = \pi^{-z/2}\Gamma(z/2)\zeta(z)$.

Let  $z=x+i y$. Then 
the function \eqref{xi} can be written as
\beq
\xi(z,\chi) = A(x,y,\chi)\exp\left\{ i \theta(x,y,\chi)  \right\},
\eeq
where
\begin{align}
A(x,y,\chi) & =
\(\dfrac{k}{\pi}\)^{\tfrac{x+a}{2}}
\left|\Gamma\(\dfrac{x+a+iy}{2}\)\right| 
\left|L(x+iy,\chi)\right|, \label{Achi} \\
\theta(x,y,\chi) &= \arg \Gamma\(\dfrac{x+a+iy}{2}\) - 
\dfrac{y}{2}\log \(\dfrac{\pi}{k}\) 
- \dfrac{1}{2}\arg G(\chi) + 
\arg L(x+iy,\chi) + \dfrac{\pi a}{4}. \label{thetachi}
\end{align}
From \eqref{xiconju} we also conclude that
$A(x,y,\chi) = A(x,-y,\chi^*)$ and $\theta(x,y,\chi) = -\theta(x,-y,\chi^*)$.
Denoting
\beq
\xi^*(1-z, \chi) = A'(x,y,\chi) \exp\left\{ -i \theta'(x,y,\chi)\right\}
\eeq
we therefore have 
\beq\label{thetaprime}
A'(x,y,\chi) = A(1-x,y, \chi), \qquad \theta'(x,y, \chi) = \theta(1-x,y,\chi).
\eeq
Taking the modulus of \eqref{FExi} we also have that
$A(x,y,\chi)=A'(x,y,\chi)$ for any $z$.

On the critical strip, the functions $L(z,\chi)$ and $\xi(z,\chi)$ 
have the same zeros.
Thus the zeros can be defined by $\xi(\rho,\chi) = 0$. Thus on a zero
we clearly have
\beq\label{onzero}
\lim_{\delta\to0^+} 
\left\{ \xi (\rho+\delta, \chi)  +  \xi^* (1-\rho-\delta,\chi ) \right\} =0.
\eeq
Let us denote
\beq
\label{Bdef}
B(x,y,\chi) = e^{i\theta (x,y,\chi)} + e^{- i\theta' (x,y,\chi)} .
\eeq
Since $A=A'$ everywhere,  from \eqref{onzero} we 
conclude that on a zero we have
\beq\label{theta_zero}
\lim_{\delta \to 0^+}  A(x+ \delta ,y, \chi) B(x+\delta ,y, \chi) = 0.
\eeq
This equation is satisfied if $A=0$, $B=0$ or both,
where the $\delta\to 0^+$ limit is implicit 
here 
and in the following.      
The function $A$ has the same zeros as $L(z,\chi)$ 
since $A(x,y,\chi) \propto |L(z,\chi)|$.    
There is more structure in $B$,  so 
let us consider its zeros.
This equation has the general solution $\theta + \theta' = (2n+1)\pi$,
which is a family of curves $y(x)$ and thus cannot correspond to the 
zeros of a complex analytic function which must be isolated points. Therefore,
let us choose the particular solution within this general class 
given by
\beq\label{particular}
\theta = \theta' , \qquad \lim_{\delta\to 0^{+}}\cos\theta = 0. 
\eeq
A more lengthy discussion of why the above particular solution 
corresponds to the Riemann zeros was given in \cite{GuiLeClair}.

Let us define the function
\beq\label{RSgen}
\begin{split}
\vartheta_{k,a} (y) &\equiv  
\arg \Gamma \( \dfrac{1}{4} + \dfrac{a}{2} + i \, \dfrac{y}{2}  \)  
- \dfrac{y}{2}  \log \( \dfrac{\pi}{k} \) \\
& = \Im\[\log\Gamma\( \dfrac{1}{4} + \dfrac{a}{2} +i \,  \dfrac{y}{2} \)\]
- \dfrac{y}{2} \log \( \dfrac{\pi}{k} \).
\end{split}
\eeq
When $k=1$ and $a=0$, the function \eqref{RSgen} is just the usual 
Riemann-Siegel $\vartheta$ function. 
Since the function $\log\Gamma$ has a complicated branch cut, one can
use the following series representation in \eqref{RSgen} \cite{Boros}
\beq
\log\Gamma(z) = -\gamma z - \log z -
\sum_{n=1}^{\infty}\left\{ \log\(1+\dfrac{z}{n}\) -\dfrac{z}{n}\right\},
\eeq
where $\gamma$ is the Euler-Mascheroni constant. Nevertheless,
most numerical packages already have the $\log\Gamma$ function implemented.

On the \emph{critical line} the first equation in \eqref{particular} is
already satisfied.  From the second equation we have 
$\lim_{\delta\to0^+}\theta=\(n+\tfrac{1}{2}\)\pi$,  therefore
\beq\label{almost_final}
\vartheta_{k,a}(y_n) + 
\lim_{\delta\to 0^{+}}\arg L\(\tfrac{1}{2}+\delta+iy_n,\chi\)
- \dfrac{1}{2}\arg G\(\chi\) + \dfrac{\pi a }{4} = 
\(n+\dfrac{1}{2}\)\pi .
\eeq
Analyzing the left hand side of \eqref{almost_final} we can see that it
has a minimum, thus we shift $n \to n - (n_0 + 1)$ for a specific $n_0$, 
to label the zeros according to the convention that the first positive 
zero is labelled by $n=1$. 
Thus the upper half of the critical line will have the zeros
labelled by $n=1,2,\dotsc$ corresponding to positive $y_n$, while the 
lower half will have the negative values $y_n$ 
labelled by $n=0,-1,\dotsc$. The integer $n_0$ depends 
on $k$, $a$ and $\chi$, and should be chosen according to each specific case.
In the cases we analyze below $n_0=0$, whereas for the trivial 
character $n_0 =1$. In practice, the value of $n_0$ can always be 
determined by  plotting \eqref{almost_final} with $n=1$, 
passing all terms to its left hand side. Then it is trivial to adjust the 
integer $n_0$ such that the graph passes through the point 
$(y_1, 0)$ for the first jump, corresponding to the first positive solution.
Henceforth we will \emph{omit} the integer $n_0$ in the equations, 
since all cases analyzed in this paper have $n_0=0$. Nevertheless, the reader
should bear in mind that for other cases, it may  be necessary
to replace $n \to n - n_0$ in the following equations.

In summary,  these zeros have the form $\rho_n = \tfrac{1}{2}+i y_n$, 
where for a given  
$n \in \mathbb{Z}$, the imaginary part $y_n$ is the  solution of 
the equation
\beq\label{exact}
\vartheta_{k,a}(y_n) + 
\lim_{\delta\to0^{+}} \arg L\(\tfrac{1}{2}+\delta+iy_n,\chi\)
- \dfrac{1}{2}\arg G\(\chi\) = \(n - \dfrac{1}{2} - \dfrac{a}{4} \)\pi.
\eeq
 It is important to note
that the above limit is defined with a positive $\delta$.   
This limit is well defined,  is generally not equal to zero,  and
consistent with other definitions of $\arg L$.    
It also controls  its wild oscillation when solving the equation numerically.

\subsection{An asymptotic equation for the \boldmath{$n$}-th zero}
\label{sec:approx}

\def\sig{\sigma} 

From Stirling's formula we have the following 
asymptotic form for $y \to \pm \infty$:
\beq
\label{approxRS}
\vartheta_{k,a} (y) = \mathrm{sgn}(y)
\left\{ \dfrac{|y|}{2} \log\(  \dfrac{k |y|}{2 \pi e} \)  +  
\dfrac{2a-1}{8} \pi
+O(1/y) \right\}.
\eeq
The first order approximation of \eqref{exact}, i.e.
neglecting terms of $O(1/y)$, is therefore given by
\begin{multline}\label{asymptotic}
\sig_n  \dfrac{|y_n|}{2\pi}\log\(\dfrac{k \, |y_n|}{2\pi e }\)
+\dfrac{1}{\pi}\lim_{\delta\to 0^{+}}
\arg L\(\tfrac{1}{2}+ \delta + i \sig_n |y_n|, \chi\)
- \dfrac{1}{2\pi}\arg G\(\chi\) \\
 = n + \dfrac{\sig_n -4 - 2a(1+\sig_n)}{8},
\end{multline}
where $\sig_n = 1$ if $n>0$ and $\sig_n = -1$ if $n\le 0$. For $n>0$
we have $y_n = |y_n|$ and for $n\le 0$ $y_n = - |y_n|$.

\subsection{An explicit approximate solution in terms of the Lambert function}
\label{sec:lambert}

Using the definition of the Lambert $W$
function, $W(z) e^{W(z)} = z$, if we neglect the much smaller  
$\arg L$ term in 
\eqref{asymptotic} we can find an exact solution. Let
\beq
A_{n}\(\chi\) = \sig_n \(n + \dfrac{1}{2\pi}\arg G(\chi)\)
+\dfrac{1 - 4\sig_n - 2a\(\sig_n + 1\)}{8}.
\eeq
Considering the transformation $|y_n| = 2 \pi A_{n} x_n^{-1}$,
equation \eqref{asymptotic} can be
written as $x_n e^{x_n} = k A_{n} e^{-1}$. Thus the approximate
solution that takes into account only the smooth part of \eqref{exact} 
is explicitly given by
\beq\label{approx}
\tilde{y}_n = \dfrac{2\pi \sig_n A_{n}\(\chi\)}{
W\[ k \, e^{-1} A_{n}\(\chi\) \]},
\eeq
where $W$ is the principal branch of the Lambert $W$ function over real
values. The $W$ function is implemented in most numerical packages, thus
\eqref{approx} can easily estimate arbitrarily  high zeros on the line.
In \eqref{approx} $n=1,2,\dotsc$ correspond to positive $y_n$
solutions,  while $n=0,-1,\dotsc$ correspond to negative $y_n$ solutions.

\subsection{Counting formulas}
\label{sec:counting}

Let us define $N^+_0 (T,\chi)$ as  the number of zeros on the critical 
line with $0< \Im( \rho) < T$ and  $N^-_0 (T,\chi)$ as the number of zeros 
with $-T < \Im( \rho)  < 0$. As explained in section \ref{sec:character},
$N^+_0(T,\chi) \neq N^-_0 (T,\chi)$
if the characters are complex numbers, since  the zeros are not 
symmetrically distributed between the upper and lower half of 
the critical line. 

\begin{figure}
\centering
\includegraphics[width=0.6\linewidth]{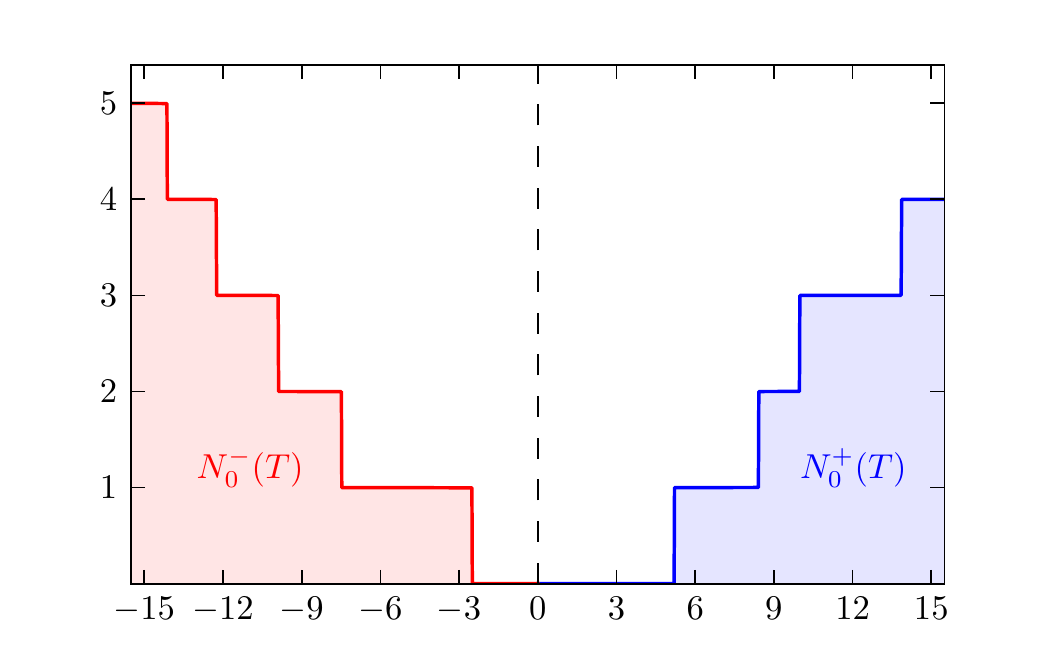}
\caption{(a) Exact counting formulae \eqref{count_exact1}  
and \eqref{count_exact2}.
Note that they are not symmetric with respect to the origin, since the
$L$-zeros for complex $\chi$ are not complex conjugates.
We used $\chi=\chi_{7,2}$ \eqref{char72}.}
\label{fig:counting}
\end{figure}

The counting formula $N^+_0(T,\chi)$  is obtained 
from \eqref{exact} by replacing $y_n \to T$ and $n \to N^+_0 + 1/2 $, 
therefore
\beq\label{count_exact1}
N^+_0(T,\chi) = \dfrac{1}{\pi}\vartheta_{k,a}(T) + 
\dfrac{1}{\pi}\arg L\(\tfrac{1}{2}+ i T, \chi\)
- \dfrac{1}{2\pi}\arg G\(\chi\) + \dfrac{a}{4}.
\eeq
Comparing with the counting formulas,  which are staircase functions,   
the left hand side of \eqref{exact} is a monotonically  
increasing function. Assuming that such a function is continuous with
the $\delta$ limit, then equation \eqref{exact} should have a unique
solution for every $n$. This justifies the passage from \eqref{exact} to
\eqref{count_exact1}.
Analogously, the counting formula on the lower half line is
given by
\beq\label{count_exact2}
N_0^-(T,\chi) = \dfrac{1}{\pi}\vartheta_{k,a}(T) - 
\dfrac{1}{\pi}\arg L \(\tfrac{1}{2} - i T, \chi\)
+\dfrac{1}{2\pi}\arg G(\chi) - \dfrac{a}{4}.
\eeq
Note that in \eqref{count_exact1} and \eqref{count_exact2} $T$ is positive.
Both cases  are plotted in \fref{fig:counting} for the character 
$\chi_{7,2}$ shown in \eqref{char72}.
One can notice that they are precisely staircase functions, jumping by one
at each zero. Note also that the functions
are not symmetric about  the origin, since for a complex $\chi$ the
zeros on upper and lower half lines are not simply complex conjugates.

From \eqref{approxRS} we also have the first order approximation
for $T \to \infty$,
\beq\label{count_approx}
N^+_0(T,\chi) = \dfrac{T}{2\pi}\log\(\dfrac{k \, T}{2\pi e }\) + 
\dfrac{1}{\pi}\arg L\(\tfrac{1}{2}+i T, \chi\)
-\dfrac{1}{2\pi}\arg G\(\chi\) - \dfrac{1}{8} + \dfrac{a}{2} .
\eeq
Analogously, for the lower half line we  have
\beq\label{count_approx2}
N^-_0(T,\chi) = \dfrac{T}{2\pi}\log\(\dfrac{k \, T}{2\pi e }\) - 
\dfrac{1}{\pi}\arg L\(\tfrac{1}{2}-i T, \chi\)
+\dfrac{1}{2\pi}\arg G\(\chi\) - \dfrac{1}{8}.
\eeq
As  in \eqref{exact}, again we are omitting $n_0$ since in the cases 
below $n_0=0$, but for other cases one may need to include $\pm n_0$ on
the right hand side of $N_0^{\pm}$, respectively.

It is known that the number of zeros on the \emph{whole critical strip} 
up to height $T$, i.e. 
$0 < x < 1$ and $0 < y < T$, is given by \cite{Montgomery}
\beq
\label{count_strip}
N^+(T,\chi) = \dfrac{1}{\pi} \vartheta_{k,a}\(T\)
+\dfrac{1}{\pi}\arg L\(\tfrac{1}{2}+iT,\chi\) -
\dfrac{1}{\pi}\arg L\(\tfrac{1}{2},\chi\).
\eeq
From Stirling's approximation and noticing that 
$2a-1=-\chi(-1)$, for $T\to\infty $ we obtain the asymptotic
approximation \cite{Selberg1,Montgomery}
\beq
\label{selberg_counting}
N^+(T,\chi) = \dfrac{T}{2\pi}\log\(\dfrac{k \, T}{2\pi e} \)
+ \dfrac{1}{\pi}\arg L\(\tfrac{1}{2} + i T,\chi\) 
- \dfrac{1}{\pi}\arg L\(\tfrac{1}{2},\chi\) 
- \dfrac{\chi(-1)}{8} + O(1/T).
\eeq
Both formulas \eqref{count_strip} and \eqref{selberg_counting} are
exactly the same as \eqref{count_exact1} and \eqref{count_approx}, 
respectively. This can be seen as follows. 
From \eqref{FExi} we conclude that $\xi$ is real on the critical line.
Thus $\arg\xi\(\tfrac{1}{2}\)=0=
-\tfrac{1}{2}\arg G\(\chi\)+\arg L\(\tfrac{1}{2},\chi\)+\tfrac{\pi a}{4}$.
Then, replacing $\arg G$ in \eqref{exact} we obtain
\beq\label{exact_v2}
\vartheta_{k,a}\(y_n\) + 
\lim_{\delta\to0^+}\arg L\(\tfrac{1}{2}+\delta+iy_n, \chi\)
-\arg L\(\tfrac{1}{2},\chi\) = \(n -\tfrac{1}{2}\)\pi.
\eeq
%Note also that we have 
%$n_0 = \tfrac{1}{2\pi}\arg G(\chi) - 
%\tfrac{1}{\pi}\arg L\(\tfrac{1}{2},\chi\) -a/4$. 
Replacing $y_n\to T$ and $n\to N_0^+ + 1/2$ in 
\eqref{exact_v2} we have precisely the expression
\eqref{count_strip}, and also \eqref{selberg_counting} for $T\to \infty$.
Therefore, we conclude that $N_0^+(T,\chi) = N^+(T,\chi)$ exactly.
From \eqref{xiconju} we see that negative zeros for character
$\chi$ correspond to positive zeros for character $\chi^*$.
Then for $-T < \Im(\rho) < 0$ the counting on the strip
also coincides  with the counting on the line, since
$N_0^-(T,\chi) = N_0^+(T,\chi^*)$ and $N^-(T,\chi) = N^+(T,\chi^*)$.
Therefore, the number of zeros on the whole \emph{critical strip}
is the same as the number of zeros on the \emph{critical line} obtained
as solutions of \eqref{exact}.

\subsection{The generalized Riemann hypothesis}
\label{sec:GRH}

Our previous conclusions were based on the particular solution
\eqref{particular} of $\lim_{\delta\to 0^+}B=0$. 
A zero of $B=0$ is also a zero of $L(z,\chi)=0$. 
The number of $L$-zeros on the  \emph{entire critical strip} is given by
\eqref{count_strip}. As discussed above,  if the equation \eqref{exact} has
a unique solution for every $n$,   
then 
$N_0^+(T,\chi) = N^+(T,\chi)$, implying that 
all zeros are  on the \emph{critical line}.
It means that \eqref{particular} captures all non-trivial zeros of $L$.
Since \eqref{exact}, or equivalently \eqref{exact_v2},  arises from 
$B=0$,  there are no additional zeros from $A=0$ 
nor from the most general solution to $B=0$ since the counting 
formula $N(T)$ is already saturated.

Another important fact follows from the one-to-one correspondence 
between the zeros of $L(z,\chi)$ with the zeros of $\cos\theta = 0$, 
which are all simple.
Thus,  under the validity of \eqref{exact} 
and \eqref{count_exact1}, the non-trivial zeros of $L(z,\chi)$ are simple.
Note that \eqref{exact} gives a different equation for each $n$, 
so there are no repeated solutions.

\subsection{Numerical solutions}
\label{sec:numerical}

We can solve the equation \eqref{exact} starting
with the approximation given by \eqref{approx}. 
We will illustrate this for some specific examples, using the root 
finder function provided by Mathematica \footnote{The Mathematica notebook 
we used to carry out these computations has only about a dozen  
lines of code and is available on the arXiv in math.NT as an ancillary file 
to this submission.}.   

The numerical procedure is carried out as follows: 
\begin{enumerate}
\item \label{step1} We solve \eqref{exact} looking for the solution in 
a region centered around the number $\tilde{y}_n$ provided 
by \eqref{approx}, with a not so small $\delta$, 
for instance $\delta \sim 10^{-5}$.
\item \label{step2} We solve \eqref{exact} again but now centered around 
the solution obtained in step \ref{step1} above, and we decrease $\delta$, for 
instance $\delta \sim 10^{-8}$.
\item We repeat the procedure in step \ref{step2} above, 
decreasing $\delta$ again.
\item Through successive iterations, and decreasing $\delta$ 
each time, it is possible to obtain solutions as accurately  as desirable.
In carrying this out,  it is important to not allow $\delta$ to be 
exactly zero.   
\end{enumerate}

We will illustrate our formulas with the primitive 
characters $\chi_{7,2}$ and $\chi_{7,3}$,   
since they possess the full  generality of $a=0$ and $a=1$ and 
complex components.  
There are actually $\varphi(7)=6$ distinct characters mod $7$.   

\paragraph*{\bf Example \boldmath{$\chi_{7,2}$}.}
Consider $k=7$ and $j=2$, i.e. we are computing
the Dirichlet character $\chi_{7,2}(n)$. For this case $a=1$.
Then we have the following components:
\beq\label{char72}
\begin{tabular}{@{}c|ccccccc@{}}
$n$             & $1$ & $2$ & $3$ & $4$ & $5$ & $6$ & $7$ \\
\midrule[0.3pt]
$\chi_{7,2}(n)$ &  
$1$ & $e^{ 2\pi i /3}$ & $e^{\pi i / 3}$ & $e^{-2\pi i / 3}$ &
$e^{-\pi i / 3}$ & $-1$ & $0$ 
\end{tabular}
\eeq
The first few zeros, positive and negative, obtained by solving
\eqref{exact} are shown in \tref{zeros_1}. The solutions shown are
easily obtained with $50$ decimal places of accuracy, and agree with
the ones in \cite{Oliveira}, which were computed up to $20$ decimal places.

\begin{table}
\def\arraystretch{0.7}
\centering
\begin{tabular}{@{}rrr@{}}
\toprule[1pt]
$n$ & $\tilde{y}_n$ & $y_n$ \\
\midrule[0.4pt]
$10$ &  $25.57$ &  $25.68439458577475868571703403827676455384372032540097$ \\
$9$ &   $23.67$ &  $24.15466453997877089700472248737944003578203821931614$ \\
$8$ &   $21.73$ &  $21.65252506979642618329545373529843196334089625358303$ \\
$7$ &   $19.73$ &  $19.65122423323359536954110529158230382437142654926200$ \\
$6$ &   $17.66$ &  $17.16141654370607042290552256158565828745960439000612$ \\
$5$ &   $15.50$ &  $15.74686940763941532761353888536874657958310887967059$ \\
$4$ &   $13.24$ &  $13.85454287448149778875634224346689375234567535103602$ \\
$3$ &   $10.81$ &   $9.97989590209139315060581291354262017420478655402522$ \\
$2$ &    $8.14$ &   $8.41361099147117759845752355454727442365106861800819$ \\
$1$ &    $4.97$ &   $5.19811619946654558608428407430395403442607551643259$ \\
$0$ &   $-3.44$ &  $-2.50937455292911971967838452268365746558148671924805$ \\
$-1$ &  $-7.04$ &  $-7.48493173971596112913314844807905530366284046079242$ \\
$-2$ &  $-9.85$ &  $-9.89354379409772210349418069925221744973779313289503$ \\
$-3$ & $-12.35$ & $-12.25742488648921665489461478678500208978360618268664$ \\
$-4$ & $-14.67$ & $-14.13507775903777080989456447454654848575048882728616$ \\
$-5$ & $-16.86$ & $-17.71409256153115895322699037454043289926793578042465$ \\
$-6$ & $-18.96$ & $-18.88909760017588073794865307957219593848843485334695$ \\
$-7$ & $-20.99$ & $-20.60481911491253262583427068994945289180639925014034$ \\
$-8$ & $-22.95$ & $-22.66635642792466587252079667063882618974425685038326$ \\
$-9$ & $-24.87$ & $-25.28550752850252321309973718800386160807733038068585$ \\
\bottomrule[1pt]
\end{tabular}
\caption{Numerical solutions of \eqref{exact} starting
with the approximation \eqref{approx}, for the
character \eqref{char72}. The solutions are
accurate to $50$ decimal places and verified to
$\left|L\(\tfrac{1}{2} + i y_n, \chi_{7,2}\)\right| \sim 10^{-50}$.}
\label{zeros_1}
\end{table}

\paragraph*{\bf Example \boldmath{$\chi_{7,3}$}.} 
Consider $k=7$ and $j=3$, such that $a=0$.
In this case the components of $\chi_{7,3}(n)$ are the following:
\beq\label{char73}
\begin{tabular}{@{}c|ccccccc@{}}
$n$             & $1$ & $2$ & $3$ & $4$ & $5$ & $6$ & $7$ \\
\midrule[0.3pt]
$\chi_{7,3}(n)$ &  
$1$ & $e^{ -2\pi i /3}$ & $e^{2 \pi i / 3}$ & $e^{2\pi i / 3}$ &
$e^{-2 \pi i / 3}$ & $1$ & $0$ 
\end{tabular}
\eeq
The first few solutions of \eqref{asymptotic} are shown in 
\tref{zeros_2} and are accurate up to $50$ decimal places, and
agree with the ones obtained in \cite{Oliveira}. 
As stated previously, the solutions to equation \eqref{exact} can be 
calculated to any desired level of accuracy. For instance, continuing with 
the character $\chi_{7,3}$, we can easily compute the following 
number for $n=1000$, accurate to $100$ decimal places, i.e. $104$ digits:
\beq\nonumber
\begin{split}
y_{1000} = 1037.&56371706920654296560046127698168717112749601359549 \\[-.8em]
                &01734503731679747841764715443496546207885576444206
\end{split}
\eeq
We also have been able to solve the equation for high zeros to 
high accuracy, up to the millionth zero, some of which are 
listed in Table \ref{highzeros},   and were previously unknown.   

\begin{table}
\def\arraystretch{0.7}
\centering
\begin{tabular}{@{}rrr@{}}
\toprule[1pt]
$n$ & $\tilde{y}_n$ & $y_n$ \\
\midrule[0.4pt]
$10$ &  $25.55$  &  $26.16994490801983565967242517629313321888238615283992$ \\
$9$ &   $23.65$  &  $23.20367246134665537826174805893362248072979160004334$ \\
$8$ &   $21.71$  &  $21.31464724410425595182027902594093075251557654412326$ \\
$7$ &   $19.71$  &  $20.03055898508203028994206564551578139558919887432101$ \\
$6$ &   $17.64$  &  $17.61605319887654241030080166645399190430725521508443$ \\
$5$ &   $15.48$  &  $15.93744820468795955688957399890407546316342953223035$ \\
$4$ &   $13.21$  &  $12.53254782268627400807230480038783642378927939761728$ \\
$3$ &   $10.79$  &  $10.73611998749339311587424153504894305046993275660967$ \\
$2$ &    $8.11$  &   $8.78555471449907536558015746317619235911936921514074$ \\
$1$ &    $4.93$  &   $4.35640162473628422727957479051551913297149929441224$ \\
$0$ &   $-5.45$  &  $-6.20123004275588129466099054628663166500168462793701$ \\
$-1$ &  $-8.53$  &  $-7.92743089809203774838798659746549239024181788857305$ \\
$-2$ & $-11.15$  & $-11.01044486207249042239362741094860371668883190429106$ \\
$-3$ & $-13.55$  & $-13.82986789986136757061236809479729216775842888684529$ \\
$-4$ & $-15.80$  & $-16.01372713415040781987211528577709085306698639304444$ \\
$-5$ & $-17.94$  & $-18.04485754217402476822077016067233558476519398664936$ \\
$-6$ & $-20.00$  & $-19.11388571948958246184820859785760690560580302023623$ \\
$-7$ & $-22.00$  & $-22.75640595577430793123629559665860790727892846161121$ \\
$-8$ & $-23.94$  & $-23.95593843516797851393076448042024914372113079309104$ \\
$-9$ & $-25.83$  & $-25.72310440610835748550521669187512401719774475488087$ \\
\bottomrule[1pt]
\end{tabular}
\caption{Numerical solutions of \eqref{exact} starting
with the approximation \eqref{approx}, for the character \eqref{char73}. 
The solutions are
accurate to $50$ decimal places and verified to
$\left|L\(\tfrac{1}{2} + i y_n, \chi_{7,3}\)\right| \sim 10^{-50}$.}
\label{zeros_2}
\end{table}

\begin{table}
\def\arraystretch{0.7}
\centering
\begin{tabular}{@{}crr@{}}
\toprule[1pt]
$n$ & $\tilde{y}_n$ & $y_n$ \\
\midrule[0.4pt]
$10^3$ & $1037.61$ & 
  $1037.563717069206542965600461276981687171127496013595490$ \\
$10^4$ & $7787.18$ &      
  $7787.337916840954922060149425635486826208937584171726906$ \\
$10^5$ & $61951.04$ &  
 $61950.779420880674657842482173403370835983852937763461400$ \\
$10^6$ & $512684.78$ & 
$512684.856698029779109684519709321053301710419463624401290$ \\
\bottomrule[1pt]
\end{tabular}
\caption{
Higher zeros for the Dirichlet character \eqref{char73}.  
These solutions to \eqref{exact}  are
accurate to $50$ decimal places.}
\label{highzeros}
\end{table}

\section{\boldmath{$L$}-functions based on modular forms}
\label{sec:modular}

\def\shalf{\tfrac{1}{2}}
\def\SL{SL_2\(\mathbb{Z}\)}
\def\xihat{\hat{\xi}}
\def\Lhat{\hat{L}}

\subsection{General level one forms}  

The \emph{modular group} can be represented by the set
of $2 \times 2$ integer matrices,
\beq \label{SL2Z}
\SL =  \left\{ 
A = \begin{pmatrix}[0.7]
a & b \\
c & d
\end{pmatrix} \, \bigg\vert \, \,  a,b,c,d \in \mathbb{Z}, \ \det A = 1
\right\},
\eeq
provided each matrix $A$ is identified with $-A$, i.e. $\pm A$ are regarded
as the same transformation. Thus for $\tau$ in the upper half 
complex plane, it transforms as 
$\tau  \mapsto A \tau = \tfrac{a \tau + b}{c\tau + d}$ 
under the action of the modular group.
A \emph{modular form} $f$ of weight $k$ is a function that is analytic in the 
upper half complex plane which satisfies the functional relation 
\cite{Apostolmodular} 
\beq \label{fform}
f\(\dfrac{a\tau  + b}{c\tau + d} \) =  \(c\tau +d\)^k   f(\tau).
\eeq
If the above equation is satisfied for all of $\SL$,  
then $f$ is referred to as being of level one.
It is possible to define higher level modular forms  which satisfy 
the above equation for a subgroup of $\SL$.   Since our results are easily
generalized to the higher level case,  henceforth we will  only consider 
level 1 forms.     

For the $\SL$ element 
$\renewcommand\arraycolsep{2pt}
\begin{pmatrix}[0.6] 1 & 1 \\ 0 & 1 \end{pmatrix}$, 
the above implies  the periodicity $f(\tau ) = f(\tau+1)$, 
thus it has a Fourier series
\beq \label{Fourier}
f(\tau ) = \sum_{n=0}^\infty  a_f (n) \, q^n  ,  
\qquad q \equiv  e^{2 \pi i \tau}.
\eeq
If $a_f(0) = 0$ then $f$ is called a \emph{cusp form}. 

From the Fourier coefficients, one can define the Dirichlet series
\beq \label{Lmod}
L_f \(z\) = \sum_{n=1}^\infty \dfrac{a_f \(n\)}{n^z}.
\eeq
The functional equation for $L_f\(z\)$ relates it to 
$L_f\(k-z\)$,  so that the critical line is $\Re (z) = \tfrac{k}{2}$, 
where $k \ge 4$ is an even integer.   One can always shift the critical line
to $\tfrac{1}{2}$ by  replacing $a_f (n)$  by  $a_f (n)/n^{(k-1)/2}$,  
however we will not do this here.   
Let us define 
\beq \label{ximod}
\Lambda_f(z) \equiv \(2\pi\)^{-z} \, \Gamma \( z \) \, L_f (z).
\eeq
Then the functional equation is given by \cite{Apostolmodular}
\beq \label{LFE}
\Lambda_f(z) = (-1)^{k/2} \Lambda_f(k-z).
\eeq

There are only two cases to consider since $\tfrac{k}{2}$ can be
an even or an odd  integer.  As in \eqref{xi} we can absorb
the extra minus sign  factor for the odd case. 
Thus we define $\xi_f(z) \equiv \Lambda_f(z)$ for
$\tfrac{k}{2}$ even, and  we have $\xi_f(z) = \xi_f(k-z)$,  
and $\xi_f(z) \equiv e^{-i \pi/2}\Lambda_f(z)$ for $\tfrac{k}{2}$ odd, 
implying $\xi_f(z) = \xi_f^*(k-z)$. 
Representing $\xi_f(z) = |\xi_f(z)| \, e^{i \theta (x,y)}$ where 
$z=x+i y$, we follow exactly the same steps as in the last section.
From the particular solution \eqref{particular} we conclude
that there are  infinite zeros
on the critical line $\Re(\rho) = \tfrac{k}{2}$ determined by
$\lim_{\delta\to0^+}\theta\(\tfrac{k}{2}+\delta,y, \chi\)=
\(n-\tfrac{1}{2}\)\pi$.
Therefore, these zeros have
the form $\rho_n = \tfrac{k}{2} + iy_n$, where $y_n$ is the solution
of the equation
\beq \label{exactmod}
\vartheta_k (y_n) + \lim_{\delta \to 0^+} 
\arg L_f \( \tfrac{k}{2} + \delta + i y_n \) =  
\( n - \dfrac{1 + (-1)^{k/2}}{4} \) \pi,
\eeq
where we have defined
\beq \label{RSL}
\vartheta_k (y)  \equiv \arg \Gamma\( \tfrac{k}{2} + i y \) - y \log 2\pi.
\eeq
Repeating the argument of the last section,   
if equation \eqref{exactmod} has a unique solution 
for every $n$,   then this implies that the number of zeros with 
imaginary part less 
that $T$ is given by
\beq \label{NTmod}
N_0\(T\) =  \inv{\pi} \vartheta_k (T) + 
\inv{\pi} \arg L_f \( \tfrac{ k}{2} + iT \) - \dfrac{1 - (-1)^{k/2}}{4}.
\eeq

In the limit of large $y_n$, neglecting terms of $O(1/y)$,
the equation \eqref{exactmod} becomes 
\beq \label{asymod}
y_n \log \(  \frac{y_n}{2 \pi e} \)  + 
\lim_{\delta \to 0^+}  \arg L_f \( \tfrac{k}{2} +\delta + i y_n\) 
= \( n - \dfrac{k+ (-1)^{k/2}}{4} \) \pi.
\eeq
If one ignores the small $\arg L_f$ term, then the approximate solution 
is given by
\beq \label{yLambertmod}
\tilde{y}_n = \dfrac{A_n \pi }{W\[ (2e)^{-1}A_n \]}, \qquad
A_n = n - \dfrac{k+(-1)^{k/2}}{4}. 
%\tilde{y}_n = \dfrac{\(n - \tfrac{k+(-1)^{k/2}}{4} \)\pi }{
%W\[ (2e)^{-1} \(n-\tfrac{k+(-1)^{k/2}}{4}\) \]}.
\eeq

\subsection{An example with weight \boldmath{$k=12$}}
\label{sec:ramanujan}

The simplest example is based on the Dedekind $\eta$-function
\beq \label{eta}
\eta (\tau) = q^{1/24} \, \prod_{n=1}^\infty (1-q^n) .
\eeq
Up to a simple factor,  $\eta$ is the inverse of the chiral  partition 
function of the free boson conformal field theory \cite{CFT},  where $\tau$ is
the modular parameter of the torus.   
The modular discriminant 
\beq \label{Delta} 
\Delta (\tau ) = \eta (\tau )^{24}  =  \sum_{n=1}^\infty \, \tau(n) \, q^n
\eeq
is a weight $k=12$ modular form.   
It is closely related to the inverse of the partition function 
of the bosonic string in $26$ dimensions,  where $24$ is the number 
of light-cone degrees of freedom \cite{String}.     
The Fourier coefficients  $\tau (n)$  correspond to  the 
Ramanujan $\tau$-function, and the first few are 
\beq
\begin{tabular}{@{}c|ccccccccc@{}}
$n$             & $1$ & $2$ & $3$ & $4$ & $5$ & $6$ & $7$ & $8$ & $9$  \\
\midrule[0.3pt]
$\tau(n)$ &  $1$ & $-24$ & $252$ & $-1472$ & $4830$ & $-6048$ & 
$-16744$ & $84480$ & $-113643$
\end{tabular}
\eeq
We then define the Dirichlet series
\beq \label{LRam}
L_{\Delta}(z)  =  \sum_{n=1}^\infty  \,  \frac{\tau (n) }{n^z}.
\eeq

\begin{figure}
\centering
\includegraphics[width=0.6\linewidth]{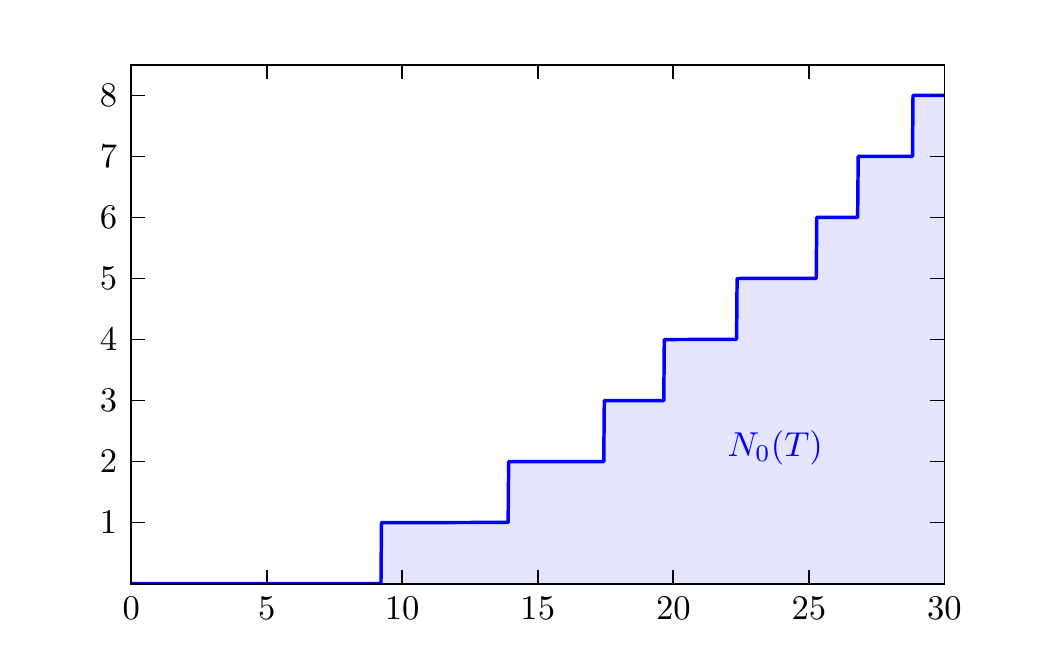}
\caption{Exact counting formula \eqref{countRam} based on 
Ramanujan $\tau$-function.}
\label{fig:ram_counting}
\end{figure}

From \eqref{exactmod} the zeros are $\rho_n = 6 + i y_n$,  
where  the $y_n$ satisfy the exact equation
\beq \label{exactRam}
\vartheta_{12} (y) + 
\lim_{\delta \to 0^{+}} \arg L_{\Delta}(6 + \delta + i y_n ) = 
\(n-\tfrac{1}{2}\) \pi.
\eeq
The counting function \eqref{NTmod} and its asymptotic approximation are
\begin{align} 
\label{countRam}
N_0(T) &= \inv{\pi} \vartheta_{12} (T) + 
\inv{\pi} \arg L_{\Delta} (6 + i T)  \\
&\approx
\dfrac{T}{\pi} \log \( \dfrac{T}{2 \pi e} \)  + \inv{\pi}\arg L_\Delta(6+iT) +
\dfrac{11}{4} + O(1/T).
\end{align}
A plot of \eqref{countRam} is shown in 
\fref{fig:ram_counting}, and we can see that it is a perfect staircase
function.

The approximate solution \eqref{yLambertmod} in terms of the 
Lambert function is given by
\beq \label{LambertRam} 
\tilde{y}_n = \dfrac{ \(n- \tfrac{13}{4} \) \pi }{W\[(2e)^{-1}
\(n-\tfrac{13}{4}\)\]}
\qquad (n=2,3,\dotsc).
\eeq
Note that the above equation is valid for $n > 1$, since $W(x)$ is not defined
for $x < -1/e$.

We follow exactly the same procedure discussed in the
beginning of section \ref{sec:numerical} to solve
equation \eqref{exactRam} numerically, starting with the approximation
provided by \eqref{LambertRam}. Some of these solutions are shown
in \tref{zerosRam} and are accurate to $50$ decimal places \footnote{See
also the Mathematica notebooks attached to this submission on arXiv.}.

\begin{table}
\def\arraystretch{0.7}
\centering
\begin{tabular}{@{}rrr@{}}
\toprule[1pt]
$n$ & $\tilde{y}_n$ & $y_n$ \\
\midrule[0.4pt]
$1$   &          &   $9.22237939992110252224376719274347813552877062243201$ \\
$2$   &  $12.46$ &  $13.90754986139213440644668132877021949175755235351449$ \\
$3$   &  $16.27$ &  $17.44277697823447331355152513712726271870886652427527$ \\
$4$   &  $19.30$ &  $19.65651314195496100012728175632130280161555091200324$ \\
$5$   &  $21.94$ &  $22.33610363720986727568267445923624619245504695246527$ \\
$6$   &  $24.35$ &  $25.27463654811236535674532419313346311859592673122941$ \\
$7$   &  $26.60$ &  $26.80439115835040303257574923358456474715296800497933$ \\
$8$   &  $28.72$ &  $28.83168262418687544502196191298438972569093668609124$ \\
$9$   &  $30.74$ &  $31.17820949836025906449218889077405585464551198966267$ \\
$10$  &  $32.68$ &  $32.77487538223120744183045567331198999909916163721260$ \\
$100$ & $143.03$ & $143.08355526347845507373979776964664120256210342087127$ \\
$200$ & $235.55$ & $235.74710143999213667703807130733621035921210614210694$ \\
$300$ & $318.61$ & $318.36169446742310747533323741641236307865855919162340$ \\
\bottomrule[1pt]
\end{tabular}
\caption{Non-trivial zeros of the modular $L$-function based
on the Ramanujan $\tau$-function, obtained from \eqref{exactRam} starting
with the approximation \eqref{LambertRam}.  
These solutions are accurate to $50$ decimal places.}
\label{zerosRam}
\end{table}

\section{Conclusion}
\label{sec:conclusion}

We have generalized the approach proposed in \cite{AL, GuiLeClair} for
the Riemann $\zeta$-function, to all primitive Dirichlet $L$-functions.
We showed that there are infinite zeros on the critical line in 
one-to-one correspondence with the zeros of $\cos\theta=0$, equation
\eqref{particular}. In this way, the zeros are enumerated and their imaginary
parts satisfy the equation \eqref{exact}.

Under the weak assumption that \eqref{exact} is well defined
for every $y_n$,  noting  that it is a monotonically  increasing function, then
it has a unique solution for every $n$. Thus one can obtain
the counting formula \eqref{count_exact1} on the \emph{critical line}. 
This agrees with the counting formula for
the number of zeros on the  \emph{entire critical strip}. Thus the zeros from
\eqref{exact} already saturate the counting formula on the strip, indicating
that all zeros must be on the critical line.

We have computed some numerical solutions of \eqref{exact} to high
accuracy, shown in section \ref{sec:numerical}. 
Thus this gives additional strong
numerical support for the validity of \eqref{exact}.

Furthermore, we have also employed the same analysis for $L$-functions
based on level one modular forms in section \ref{sec:modular}, 
and considered a specific example based on the
Ramanujan $\tau$-function. In appendix \ref{sec:general} we
unify this analysis for a more general class of $L$-functions.

\section*{Acknowledgments}
We wish to thank Tim Healey and Wladyslaw Narkiewicz for useful 
discussions.
GF is grateful for the support of CNPq  under the 
``Ci\^ encias sem fronteiras''
program.
This work is supported by the National Science
Foundation of the United States of America under grant number NSF-PHY-0757868.

%\bigskip

\appendix

\section{General \boldmath{$L$-functions}}
\label{sec:general}

With the aim  of unifying the previous results,
we extend  the main equations to a more general,  but not specific, 
class of $L$-functions.     We must of course assume that the analysis  
in section IIB,     and also discussed more throughly 
in \cite{GuiLeClair},  is valid,  and this must be checked case by case.
Thus one must bear in mind that  since we are not specific about 
the $L$-function,
but just assume some elementary properties,   there is no guarantee that 
the latter analysis is valid for every $L$-function with the properties below.  
We are thus simply  going to assume that the analysis in section 
\ref{sec:equation}  is valid 
and present the resulting equations that  would follow.  

We are going to consider $L$-functions with the properties
outlined in \cite[Chap. 5]{Iwaniec2}.  
The general $L$-function has a Dirichlet series
\beq
L(z, f) = \sum_{n=1}^{\infty}\dfrac{\lambda_f(n)}{n^{z}},
\eeq
where $\lambda_f(1) = 1$ and $\lambda_f(n)$ is a complex number. This series
is convergent for $\Re(z) > 1$, it has an Euler product of
degree $d \ge 1$, where $d$ is an integer, and admits an analytic continuation
in the whole complex plane, except for poles at $z=0$ and $z=1$. 
The arithmetic object $f$ has no specific meaning here.
For instance, $f$ can be a modular or cusp form, or it can be associated
to Dirichlet characters $\lambda_f(n) = \chi(n)$.
The object $f$ defines the particular class of $L$-functions. Let
\beq
\gamma\( z, f\) = \pi^{-d z / 2} \prod_{j=1}^{d}\Gamma\( \dfrac{z + a_j}{2} \),
\eeq
where the complex numbers $a_j$, which come in conjugate pairs, are the
so called local parameters at infinity. Let us define
\beq
\Lambda\(z, f\) \equiv k^{z/2} \gamma(z,f) \, L(z,f),
\eeq
where $k=k(f) \ge 1$ is an integer, the \emph{conductor} of $L(z,f)$. Then
it satisfies the functional equation \cite{Iwaniec2}
\beq
\Lambda(z,f) = \alpha(f) \Lambda^*(1-z,f) 
\equiv \alpha(f) \( \Lambda(1-z^*,f) \)^{*}.
\eeq
Here $\alpha$ is a 
complex phase, i.e. $|\alpha| = 1$.
The symbol $f^*$ denotes the dual of $f$, associated to the
Dirichlet series with $\lambda_{f^*}(n) = \(\lambda_f(n)\)^*$. We also
have the relations $\gamma(z,f^*) = \gamma(z, f)$ and $k(f^*)=k(f)$.
%and $\alpha(f^*)=\alpha(f)$.

We can
write the functional equation in a more symmetric form by introducing
\beq
\xi(z,f) \equiv e^{-i\beta/2}\Lambda(z,f),
\eeq
where $\alpha(f) = e^{i\beta(f)}$. Then we have the functional equation
\beq
\xi(z,f) = \xi^*(1-z,f) = \(\xi(1-z^*,f)\)^*.
\eeq
Writing the polar form $\xi(z,f) = A(x,y,f)e^{i\theta(x,y,f)}$ we have
\begin{align}
A(x,y,f) &= k^{x/2}\pi^{-dx/2} \left| L(x+iy, f)\right| 
\prod_{j=1}^{d}\left| \Gamma\( \dfrac{x+a_j + iy}{2} \) \right|, \\
\theta(x,y,f) &= \sum_{j=1}^d \arg\Gamma\( \dfrac{x+a_j + i y}{2} \)
+ \arg L(x+iy,f) + \dfrac{y}{2}\log\(\dfrac{k}{\pi^d}\) - \dfrac{\beta}{2}.
\end{align}
Denoting $\xi^*(1-z,f)=A'(x,y,f)e^{-i\theta'(x,y,f)}$, we then have
$A'(x,y,f) = A(1-x,y,f)$ and $\theta'(x,y,f)=\theta(1-x,y,f)$.
Let us define the generalized Riemann-Siegel $\vartheta_{k, a_j}$ function
\beq\label{RS_gen}
\vartheta_{k,a_j}(y) \equiv \arg\Gamma\(\dfrac{1}{4} + \dfrac{a_j}{2} + 
i \dfrac{y}{2}\) - \dfrac{y}{2}\log\(\dfrac{\pi}{k^{1/d}}\).
\eeq
Following the same previous analysis, by imposing the identity
\beq
\lim_{\delta \to 0^+}\( \xi(\rho + \delta,f) + \xi^*(1-\rho-\delta,f) \) = 0
\eeq
where $\rho = x+iy$ is a non-trivial $L$-zero, we take the particular
solution $\theta = \theta'$ and $\lim_{\delta\to0^+}\cos\theta = 0$. 
Thus $\theta = \theta'$ is satisfied by $\Re(\rho) = \tfrac{1}{2}$ and
then $\lim_{\delta\to0^+}\cos\theta\(\tfrac{1}{2},y,f\) = 0$ yields the 
equation for the $n$-th zero.
Introducing a shift $ n \to n - (n_0+1)$, where $n_0$ should 
be determined by each specific case according to the convention that
the first positive zero is labelled by $n=1$ (we will omit $n_0$ in 
the following), 
we then conclude that these zeros
have the form $\rho_n = \tfrac{1}{2} + i y_n$ 
for $n \in \mathbb{Z}$, and $y_n$ is the solution of the equation
\beq \label{exact_gen}
\sum_{j=1}^{d}\vartheta_{k,a_j}\(y_n\) + 
\lim_{\delta \to 0^+}\arg L\(\tfrac{1}{2}+\delta + iy_n, f\) - 
\dfrac{\beta}{2} = \( n 
%- n_0 
- \tfrac{1}{2} \)\pi.
\eeq
It is also possible to replace $\beta$ by noting that 
$\xi(\tfrac{1}{2}+iy, f)$ is real, thus $\arg\xi\(\tfrac{1}{2}\) = 0$ then
\beq
\dfrac{\beta(f)}{2} = \sum_{j=1}^d\arg \Gamma\( 
\dfrac{1}{4}+\dfrac{a_j}{2} \) + \arg L\(\tfrac{1}{2},f\).
\eeq
For real $a_j$ the first term vanishes. If $f^* = f$ then $L(z,f)$ is said
to be self-dual. If besides this $\alpha(f) = -1$, then 
$L\(\tfrac{1}{2}, f\) = 0$.

The counting formula on the critical line, 
for $ 0 < \Im(z) < T$, can be obtained 
from \eqref{exact_gen} by replacing $y_n \to T$ and 
$n \to N_0^+(T,f) + \tfrac{1}{2}$, thus
\beq\label{exact_counting_gen}
N_0^+(T,f) = \dfrac{1}{\pi}\sum_{j=1}^d \left\{ \vartheta_{k,a_j}(T) - 
\arg \Gamma\( \dfrac{1}{4} + \dfrac{a_j}{2} \) \right\}
+ \dfrac{1}{\pi} \arg L\(\tfrac{1}{2}+iT, f\) - 
\dfrac{1}{\pi}\arg L\(\tfrac{1}{2},f\) 
%+ n_0.
.
\eeq
The same counting on the whole strip $N^+(T,f)$ can be obtained
through the standard Cauchy's argument principle \cite{Iwaniec2}. Thus
$N_0^+(T,f) = N^+(T,f)$, justifying that the particular solution captures
all non-trivial zeros on the critical strip, therefore they must be all
on the critical line. The counting on the negative half line
$-T < \Im(z) < 0$ can be obtained from $N_0^-(T,f) = N_0^+(T,f^*)$. 
Expanding \eqref{RS_gen} from Stirling's
formula we have
\beq\label{RSasymp_gen}
\vartheta_{k,a_j}(y) = \dfrac{y}{2}\log\(\dfrac{ k^{1/d} \, y}{2\pi e }\) + 
\dfrac{2a_j-1}{8}\pi + O(1/y).
\eeq
Then from \eqref{exact_counting_gen} we have
\begin{multline}
N_0^+(T,f) = \dfrac{T}{2\pi}\log\( \dfrac{k \, T^d}{(2\pi e)^d} \)
+ \sum_{j=1}^d \left\{ \dfrac{2a_j-1}{8} - 
\inv{\pi}\arg \Gamma\( \dfrac{1}{4} +\dfrac{a_j}{2}\)  \right\} \\
+\dfrac{1}{\pi}\arg L\( \tfrac{1}{2}+iT, f \) - 
\dfrac{1}{\pi} L\( \tfrac{1}{2},f \)
%+ n_0.
+O(1/y)
.
\end{multline}

Using \eqref{RSasymp_gen} in \eqref{exact_gen} and neglecting the small
$\arg L\(\tfrac{1}{2} + iT, f\)$ term, it is possible to obtain an approximate
solution in closed form, which reads
\beq
\tilde{y}_n = \dfrac{2\pi A_n }{d \, W\[  k^{1/d} (d\, e)^{-1} A_n \]},
\eeq
where $W(x)$ is the principal value of the Lambert function over real
values, and
\beq
A_n = n - \dfrac{1}{2} %- n_0 
+ \dfrac{1}{\pi}\arg L\(\tfrac{1}{2}, f\)
+\sum_{j=1}^{d}\left\{ 
\dfrac{1}{\pi}\arg \Gamma\( \dfrac{1}{4}+\dfrac{a_j}{2} \) - 
\dfrac{2a_j-1}{8} \right\}.
\eeq

\end{document}